\documentclass[twocolumn]{svjour3}  
\smartqed  

\usepackage[dvipdfmx]{color}
\usepackage{graphicx}
\usepackage[T1]{fontenc}
\usepackage{mathptmx}  
\usepackage[scaled]{helvet}  
\usepackage{courier}
\usepackage[subrefformat=parens]{subcaption}

\makeatletter
\let\cl@chapter\undefined
\makeatletter
\usepackage[numbers,sort&compress]{natbib}
\usepackage{amsmath,amssymb}   
\usepackage{mathtools} 
\usepackage{cleveref}  
\Crefname{equation}{Eq.}{Eqs.}%
\Crefname{figure}{Fig.}{Figs.}%
\usepackage{amsbsy}
\usepackage{bm}
\usepackage{algorithmic}
\usepackage{algorithm}

\usepackage{ulem}

\usepackage{tikz}
\usetikzlibrary{calc}
\newcommand{\dif}{\mathrm{d}}
\newcommand{\fd}[2]{\frac{\dif #1}{\dif #2}}
\newcommand{\fdn}[3]{\frac{\dif^#1 #2}{\dif #3 ^#1}}
\newcommand{\lrp}[1]{\left( #1 \right)}
\newcommand{\lrs}[1]{\left[ #1 \right]}
\newcommand{\lrc}[1]{\left\{ #1 \right\}}

\begin{document}

\title{
    On the Uniqueness of the Approximate Analytical Solution of the Surf-Riding Threshold in the IMO Second Generation Intact Stability Criteria
}

\author{
    Masahiro Sakai \and Naoya Umeda
}

\institute{
    Masahiro Sakai \and Naoya Umeda \at
    Graduate School of Engineering, Osaka University, 
    2-1 Yamadaoka, Suita, Osaka 565-0871, Japan \\
    \email{sakai@naoe.eng.osaka-u.ac.jp}
}

\date{Received: date / Accepted: date}

\maketitle

\begin{abstract}
The International Maritime Organization developed the Second Generation Intact Stability Criteria (SGISC) with the objective of mitigating maritime incidents attributable to roll motion. Given that surf-riding precedes broaching, the SGISC uses vulnerability criteria for broaching based on surf-riding dynamics. The surf-riding threshold is derived by Melnikov's method in the level-2 vulnerability criterion. In this paper, the authors discuss the uniqueness of the approximate analytical solutions to be used as the surf-riding threshold in the SGISC and based on the relationship between Melnikov's analysis and its physical interpretation.
\end{abstract}

\section{Introduction}
\label{sec:intro}

The International Maritime Organization (IMO) developed the Second Generation Intact Stability Criteria (SGISC) with the objective of mitigating maritime incidents attributable to roll motion.
The interim guidelines for the SGISC were completed by the IMO's Subcommittee on Ship Design and Construction in 2020 \cite{Imo2020-bp} and are currently undergoing a trial period.
The SGISC assess five dynamic stability failure modes, including broaching, parametric roll, pure loss of stability, dead ship condition, and excessive acceleration.
Broaching occurs when a ship cannot maintain a straight course despite maximum steering efforts, with the resulting centrifugal force causing a severe heel moment.
Given that surf-riding precedes broaching, the SGISC uses level-1 and level-2 vulnerability criteria for broaching based on surf-riding dynamics \cite{Imo2020-bp,Imo2022-uf}.
For further theoretical details on the vulnerability criteria, refer to \cite{Umeda2010-xx,Maki2024-ix}.

Grim \cite{Grim1951-kk} indicated that a surf-riding boundary coincides with the trajectory connecting two unstable equilibrium points. 
This phenomenon is termed a heteroclinic bifurcation in nonlinear dynamical systems theory. 
Makov \cite{Makov1969-ya} substantiated Grim's theory by illustrating surf-riding phase planes and discovered that a self-propelled ship experiences surf-riding independently of the initial phase plane condition once the heteroclinic bifurcation occurs. 
Consequently, the heteroclinic bifurcation is regarded as the surf-riding threshold in the SGISC. 
Previous research has employed Melnikov's method \cite{Holmes1980-gy} as an approximate analytical solution method for estimating heteroclinic bifurcations. 
Kan \cite{Kan1990-lu} used Melnikov's method \cite{Holmes1980-gy} on an uncoupled surge model with locally linearised calm-water resistance. Spyrou \cite{Spyrou2001-fy} provided an exact analytical solution for the heteroclinic bifurcation of an uncoupled surge model with quadratic calm-water resistance, and additionally, Spyrou \cite{Spyrou2006-ds} applied Melnikov's method to a surge model with cubic calm-water resistance.
Maki et al. \cite{Maki2010-dy,Maki2014-rt} concluded these studies by applying Melnikov's method to an uncoupled surge model with general polynomial calm-water resistances, determining both lower and upper surf-riding thresholds, and validating these results through numerical bifurcation analysis and free-running model experiments. For further experimental validation of surf-riding estimates or the level-2 vulnerability criterion, refer to \cite{Yamamura2009-bt,Ito2014-hb}.

The IMO adopted the method proposed by Maki et al. \cite{Maki2010-dy,Maki2014-rt} for the level-2 vulnerability criterion \cite{Imo2020-bp,Imo2022-uf}, recommending the use of a numerical iteration method to solve the surf-riding threshold derived by Melnikov's method. 
However, if we use a numerical iteration method, some uncertainties remain, such as the uniqueness and convergence of the solutions. 
Sakai et al. \cite{Sakai2017-ns} remarked that the surf-riding threshold can be obtained as a solution to a quadratic equation when the propeller thrust coefficient is represented by a quadratic polynomial. Responding to their work, the approach based on the solution to a quadratic equation was added to the explanatory notes to the SGISC \cite{Imo2022-uf}. It could make us possible to exclude any uncertainties of the numerical iteration method. 
For this purpose, this study attempts to prove the uniqueness of the approximate analytical solution to be used as the surf-riding threshold in the SGISC and examines the relationship between the analytical solution and its physical interpretation.

The initial results of the investigation presented in this paper were previously described by Sakai et al. \cite{Sakai2017-ns}. In this paper, the results are presented more extensively, with more details, and also with new remarks.
    
\section{Surf-Riding Threshold by Melnikov's method}\label{sec:SREst}
In this section, the authors briefly explain the level-2 vulnerability criterion for broaching in the SGISC.
The uncoupled surge equation used is 
\begin{equation} \label{eq:me}
    \left(m+m_x\right) \ddot{\xi_\mathrm{G}}+\left[R\left(u\right)-T_\mathrm{e} \left(u ; n\right)\right]-X_\mathrm{w}\left( \xi_\mathrm{G} \right)=0
\end{equation}
Here, $m$ is the mass of the ship, 
$m_\mathrm{x}$ is the added mass in surge of the ship, 
$\xi_\mathrm{G}$ is the relative position of the ship centre of gravity from a wave trough where the wave propagating direction is set as positive, 
$R\left(u\right)$ is the calm-water resistance as a function of the ship instantaneous forward speed $u$, 
$T_\mathrm{e} \left(u ; n\right)$ is the effective propeller thrust as a function of $u$ and number of revolutions of the propeller(s) $n$, and 
$X_\mathrm{w}$ is the wave-induced surge force. 
In this equation, the variation in thrust due to wave particle velocity and higher-order terms, such as the added resistance in waves, are ignored.

When a ship surf-rides a wave at $\xi_\mathrm{G}=\xi_\mathrm{G, SR}$, the ship runs at the same speed as the wave celerity $c_\mathrm{w}$, and the resistance, effective thrust and wave-induced surge force are balanced. In other words, the existence of $n$ satisfying \Cref{eq:equilibrium} is a necessary condition for surf-riding.
\begin{equation} \label{eq:equilibrium}
    - R\left(c_\mathrm{w}\right) + T_\mathrm{e} \left(c_\mathrm{w} ; n\right) + X_\mathrm{w}\left( \xi_\mathrm{G, SR} \right)=0
\end{equation}
Note that surf-riding could occur for any initial conditions if $n$ is in a certain range, and
the boundaries are global bifurcation points, termed the heteroclinic bifurcation. 
In the level-2 vulnerability criterion, the ship is considered vulnerable to broaching if $n$ exceeds $n_\mathrm{cr}$. Here, $n_\mathrm{cr}$ represents the lower value among the revolution numbers of the propeller(s) that correspond to the two heterophilic bifurcation points.

Melnikov's method is used to obtain the heteloclinic bifurcation point, which is the surf-riding threshold.
It is an approximate analytical method for obtaining homoclinic or heteroclinic orbits in nonlinear dynamical systems. 
Maki et al. \cite{Maki2010-dy} derived the equations of $n_\mathrm{cr}$ using Melnikov's method.
They approximated $R\left(u\right)$ and $T_\mathrm{e}(u;n)$ with $N$th-degree polynomials of $u$ and $J$ as
\begin{align} 
    R(u)=&\sum_{i=0}^N r_i u^i \label{eq:poly1}\\
    J(u;n)=&\frac{\left(1-w_\mathrm{p}\right) u}{n D} \label{eq:J} \\
    K_T\left(J(u;n)\right)=&\sum_{i=0}^N \kappa_i J(u;n)^i\\
    T_\mathrm{e}(u ; n) =&\left(1-t_\mathrm{p}\right) \rho n^2 D^4 K_T\left(J(u;n)\right) \nonumber \\
     =&\sum_{i=0}^N \kappa_i {\left(1-t_\mathrm{p}\right)\left(1-w_\mathrm{p}\right)^i \rho D^{4-i}}{n^{2-i} u^i} \label{eq:poly2}
\end{align}
Here, $w_\mathrm{p}$ is the effective wake fraction, $t_\mathrm{p}$ is the thrust deduction coefficient, $\rho$ is the density of the fluid, $D$ is the propeller diameter, $J$ is the advance constant and $K_\mathrm{T}(J)$ is the propeller thrust coefficient as a function of $J$. 
Here, $u$ is calculated as
\begin{equation} \label{eq:speed}
    u = c_\mathrm{w} + \dot{\xi_\mathrm{G}}
\end{equation}
Assuming 
\begin{equation} \label{eq:wf}
    X_\mathrm{w}\simeq -f \sin \left(k \xi_\mathrm{G} \right)    
\end{equation}
and substituting \Cref{eq:poly1,eq:poly2,eq:speed,eq:wf} into \Cref{eq:me} yield
\begin{align} \label{eq:me2}
    \left(m+m_x\right) \ddot{\xi_\mathrm{G}} &
    +\sum_{i=1}^N \sum_{j=1}^i c_{i}(n) \left(\begin{array}{c}i\\j\end{array}\right) c_\mathrm{w} ^{i-j} \dot{\xi_\mathrm{G}}{}^j \nonumber \\
    & \qquad = - R\left(c_\mathrm{w}\right) + T_\mathrm{e} \left(c_\mathrm{w} ; n\right) -f \sin \left(k \xi_\mathrm{G} \right)
\end{align}
where
\begin{equation} \label{eq:ci}
    c_i(n) = -\frac{\left(1-t_\mathrm{p}\right)\left(1-w_\mathrm{p}\right)^i \rho \kappa_i}{n^{i-2} D^{i-4}}+r_i
\end{equation}
Here, $\left(\begin{array}{c}i\\j\end{array}\right) $ indicates the binomial coefficient. 
Applying Melnikov's method to \Cref{eq:me2} yields
\begin{equation} \label{eq:Melnikov}
    2 \pi \frac{T_\mathrm{e}\left(c_\mathrm{w} ; n_\mathrm{cr}\right)-R\left(c_\mathrm{w}\right)}{f}=\sum_{i=1}^N \sum_{j=1}^i C_{i j}\left(n_\mathrm{cr}\right) \cdot (-2)^j I_j
\end{equation}
where
\begin{align}
    C_{i j}\left(n_\mathrm{cr}\right)=&c_i\left(n_\mathrm{cr}\right) \cdot \frac{1}{f k^j}
    \left(\begin{array}{c}
        i \\
        j
    \end{array}\right) 
    \frac{(fk)^{\frac{j}{2}}}{\left(m+m_x\right)^{\frac{j}{2}}} c_\mathrm{w}^{i-j}\\
    I_j=&2 \sqrt{\pi} \frac{\Gamma\left(\frac{j+1}{2}\right)}{\Gamma\left(\frac{j+2}{2}\right)}
\end{align}
Here, $f$ is the amplitude of the wave-induced surge force, $k$ is the wave number calculated as $\frac{2\pi}{\lambda}$, $\lambda$ is the wave length, and $\Gamma$ indicates the Gamma function.

$R\left(u\right)$ and $K_\mathrm{T}(J)$ are approximated by quintic and quadratic polynomials, respectively, in the level-2 vulnerability criterion.
    
\section{Quadratic Equation for the surf-riding threshold} \label{sec:Quad}
$K_\mathrm{T}(J)$ is often approximated by a quadratic polynomial in the field of naval architecture, namely,
\begin{equation} \label{eq:kappa}
    \kappa_i= 0 \qquad \left( i \ne 0, 1, 2 \right)
\end{equation}
In the context of surf-riding, it is reasonable to assume $u,n \in \mathbb{R}_{\ge 0}$, and in this region, the values of $\kappa_i$ satisfy with
\begin{align}
    \kappa_0 > & 0\\
    \kappa_2 < & 0
\end{align}
$\kappa_0$ should be positive because $K_\mathrm{T}$ should be positive when $J=0$, that is, $u=0$ and $n>0$.
$\kappa_2$ should be negative because $\kappa_2$ should represent the $K_\mathrm{T}$ reduction when $J$ is large enough.
Then, \Cref{eq:poly2} can be rewritten as
\begin{equation} \label{eq:Te}
    T_\mathrm{e}(u ; n)=\tau_0 n^2+\tau_1 nu+\tau_2 u^2
\end{equation}
where
\begin{align} 
    \tau_0=&\kappa_0\left(1-t_\mathrm{p}\right) \rho D^4 & > 0 \label{eq:tau0}\\
    \tau_1=&\kappa_1\left(1-t_\mathrm{p}\right)\left(1-w_\mathrm{p}\right) \rho D^3 & {} \label{eq:tau1}\\
    \tau_2=&\kappa_2\left(1-t_\mathrm{p}\right)\left(1-w_\mathrm{p}\right)^2 \rho D^2 & < 0 \label{eq:tau2}
\end{align}
Assuming $N>2$, \Cref{eq:ci} becomes
\begin{equation} \label{eq:ci2}
    c_i(n) = \left\{ \begin{array}{cl} 
    - \tau_1 n +r_1 & \quad (i=1)\\
    - \tau_2   +r_2 & \quad (i=2)\\
    r_i & \quad (i=3,4,\ldots,N)\\
    0 & \quad (i \ne 1,2,\ldots,N)
    \end{array}\right.
\end{equation}
The right-hand side of \Cref{eq:Melnikov} can be rewritten as
\begin{align} \label{eq:right-hand}
    &\sum_{i=1}^N \sum_{j=1}^i C_{i j} \left(n_\mathrm{cr}\right) \cdot (-2)^j I_j \nonumber\\
    &=8\frac{\tau_1 n_\mathrm{cr} - r_1}{\sqrt{f k\left(m+m_x\right)}} \nonumber\\
    &+ \sum_{j=1}^2 \left( - \tau_2+r_2 \right) \frac{1}{f k^j}\left(\begin{array}{c}i\\j\end{array}\right) 
        \frac{(fk)^{\frac{j}{2}}}{\left(m+m_x\right)^{\frac{j}{2}}} c_\mathrm{w}^{i-j} \cdot (-2)^j I_j \nonumber \\
    &+  \sum_{i=3}^N \sum_{j=1}^i  r_i \frac{1}{f k^j}\left(\begin{array}{c}i\\j\end{array}\right) 
        \frac{(fk)^{\frac{j}{2}}}{\left(m+m_x\right)^{\frac{j}{2}}} c_\mathrm{w}^{i-j} \cdot (-2)^j I_j
\end{align}
As a result of assuming $K_\mathrm{T}\lrp{J}$ is quadratic,
the right-hand side of \Cref{eq:Melnikov} becomes a 1st-degree polynomial of $n_\mathrm{cr}$, and the left-hand side of \Cref{eq:Melnikov} becomes a 2nd-degree polynomial of $n_\mathrm{cr}$.
Therefore, the surf-riding threshold can be estimated using a quadratic equation.
Note that the discussion in this section is independent of the degree of polynomial of the calm-water resistance $R\left(u\right)$.
Consequently, the surf-riding threshold can be derived by solving a quadratic equation in the level-2 vulnerability criterion.
    
\section{Physical Interpretation of Melnikov's Method for Surf-Riding Threshold}
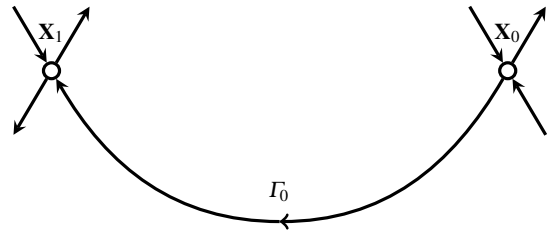
\begin{figure}[b]
    \centering
    \begin{tikzpicture}
        \node (X1) at (-3,2) {};
        \node at ($(X1) + (0, 0.5)$) {$\mathbf{X}_1$};
        \node (X0) at (3,2) {};
        \node at ($(X0) + (0, 0.5)$) {$\mathbf{X}_0$};
        \draw[very thick] (X1) circle (3pt);
        \draw[very thick] (X0) circle (3pt);
        \draw[very thick,->,>=stealth]  (3.5,1.15) -- (X0);
        \draw[very thick,->,>=stealth]  (2.5,2.85) -- (X0);
        \draw[very thick,->,>=stealth]  (X0) -- (3.5,2.85);
        \draw[very thick,->,>=stealth]  (X1) -- (-2.5,2.85);
        \draw[very thick,->,>=stealth]  (X1) -- (-3.5,1.15);
        \draw[very thick,->,>=stealth]  (-3.5,2.85) -- (X1);
        \draw[very thick, ->] (X0) to [out=-120,in=0] (0,0);
        \draw[very thick, ->,>=stealth] (0,0) to [out=180,in=-60] (X1);
        \node at (0,0.4) {$\Gamma_0$};
    \end{tikzpicture}
    \caption{Schematic view of heteroclinic orbit}
    \label{fig:schematicview}
\end{figure}

Substituting \Cref{eq:speed} into \Cref{eq:me} yeilds
\begin{equation} \label{eq:mefk}
    \left(m+m_x\right) \ddot{\xi_\mathrm{G}}+\left[R\left(u\right)-T_\mathrm{e} \left(u ; n\right)\right]+f \sin \left(k \xi_\mathrm{G} \right)=0
\end{equation}
and this can be re-written as
\begin{equation} \label{eq:mefkndb}
    \frac{1}{\sqrt{\frac{f k}{m+m_x}}^2} k \ddot{\xi_\mathrm{G}}+
    \frac{R\left(u\right)-T_\mathrm{e} \left(u ; n\right)}{f} + \sin \left(k \xi_\mathrm{G} \right)=0
\end{equation}
Putting
\begin{align} 
    y=   & k \xi_\mathrm{G}\\
    \tau=& \sqrt{\frac{f k}{m+m_x}} ~t
\end{align}
nondimensionalizes \Cref{eq:mefkndb} as
\begin{equation} \label{eq:mefknd}
   \fdn{2}{y}{\tau} + \sin y = \frac{T_\mathrm{e} \left(u ; n\right) - R\left(u\right)}{f}
\end{equation}
Here, $u$ is still used in \Cref{eq:mefknd} in order to make the following discussion clear.
The Hamilonian part of \Cref{eq:mefknd} is
\begin{equation} \label{eq:Ham}
    \fdn{2}{y}{\tau} + \sin y =0
\end{equation}
and the trajectory that connects $y=\pm \pi$ on the lower side of the vector field is
\begin{equation} \label{eq:melorbit}
    \fd{y}{\tau} = - 2 \cos \lrp{\frac{y}{2}} \triangleq \nu
\end{equation}
\Cref{fig:schematicview} shows the schematic view of the heteroclinic orbit of \Cref{eq:mefknd}.
Here, $\mathbf{X}_0, \mathbf{X}_1 \in \mathbb{R}^2$ are saddle-type equilibria (unstable surf-riding), and
$\Gamma_0$, which connects $\mathbf{X}_0$ and $\mathbf{X}_1$, is the heteroclinic orbit.
Maki et al. \cite{Maki2010-dy} obtained the Melnikov function $M\lrp{n}$ by approximating the heteroclinic orbit by the trajectory given by \Cref{eq:melorbit} as
\begin{equation} \label{eq:melfunc}
    M\lrp{n} \triangleq \int_{-\infty}^\infty \nu \lrp{\frac{T_\mathrm{e} \left(u ; n\right)-R\left(u\right)}{f}} \dif{\tau}
\end{equation}
In this integral, the integral interval represents that the heteroclinic orbit takes infinite time, and
the integrand is proportional to the power generated by $T_\mathrm{e}$ and $R$.
Under the surf-riding, the ship complies with \Cref{eq:equilibrium}.
Substituting \Cref{eq:wf} into \Cref{eq:equilibrium} yeilds
\begin{equation} \label{eq:equilibrium2}
    - R\left(c_\mathrm{w}\right) + T_\mathrm{e} \left(c_\mathrm{w} ; n\right) - f \sin \left( \xi_\mathrm{G, SR} \right)=0
\end{equation}
If $n$ is enough large to satisfy \Cref{eq:equilibrium2}, the unstable surf-riding $\mathbf{X}_l$ ($l$ is an integer) is
\begin{equation}
    \mathbf{X}_l =
    \lrp{
        \begin{array}{c}
             \xi_\mathrm{G, SR} \\
             \dot{\xi}_\mathrm{G, SR}
        \end{array}
        } 
    = \lrp{
        \begin{array}{c}
             \arcsin \lrp{ \frac{T_\mathrm{e}\lrp{ c_\mathrm{w};n} - R\lrp{c_\mathrm{w}}}{f}} -2\pi l \\
             0
        \end{array}
        }
\end{equation}
Since the potential energy due to $f \sin\lrp{k\xi_\mathrm{G}}$ at the unstable equilibria is the same, \Cref{eq:mefknd} requires that the work done by the calm-water resistance $R$ and the effective propeller thrust $T_\mathrm{e}$ must be zero in total along the heteroclinic bifurcation. Therefore, in the framework of Melnikov's method, the heteroclinic bifurcation (surf-riding threshold $n_\mathrm{cr}$) can be estimated by $M\lrp{n_\mathrm{cr}}=0$. Therefore, $n_\mathrm{cr}$ is derived with
\begin{equation} \label{eq:melbalance}
    \int_{-\infty}^\infty \nu T_\mathrm{e}\lrp{u;n_\mathrm{cr}} \dif{\tau} = \int_{-\infty}^\infty \nu R\lrp{u} \dif{\tau}
\end{equation}
Substituting \Cref{eq:melorbit} into \Cref{eq:melbalance}, changing variables, and dividing the both-sides with $-2\pi$ yield
\begin{equation} \label{eq:melbalance2}
    \frac{1}{2\pi} \int_{-\pi}^\pi T_\mathrm{e}\lrp{u;n_\mathrm{cr}} \dif{y} = \frac{1}{2\pi} \int_{-\pi}^\pi R\lrp{u} \dif{y}
\end{equation}
The ship instantaneous forward speed $u$ used for the integral in \Cref{eq:melbalance2} is defined by \Cref{eq:speed} and is transformed in the framework of Melnikov's method as
\begin{align} \label{eq:uint}
    u & = c_\mathrm{w} + \dot{\xi}_\mathrm{G} \notag \\
    &= c_\mathrm{w} + \sqrt{\frac{f}{k\lrp{m+m_\mathrm{x}}}} \fd{y}{\tau} \notag \\
    &= c_\mathrm{w} -2 \sqrt{\frac{f}{k\lrp{m+m_\mathrm{x}}}} \cos\lrp{\frac{y}{2}}
\end{align}
Based on the linear deep-water dispersion relation, $c_\mathrm{w}$ is 
\begin{equation} \label{eq:wavec}
    c_\mathrm{w} = \sqrt{\frac{g}{k}}
\end{equation}
Here, $g$ is the gravitational acceleration.

\section{Uniqueness of the surf-riding threshold} \label{sec:unique}
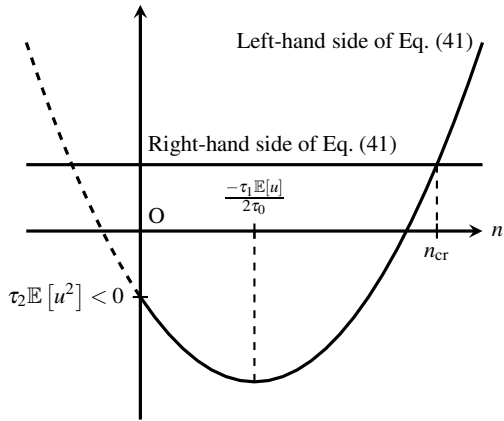
\begin{figure}[h]
    \centering
    \begin{tikzpicture}
        \draw[->,>=stealth,very thick] (-1.5,0)--(4.5,0) node[right]{$n$}; 
        \draw[->,>=stealth,very thick] (0,-2.5)--(0,3); 
        \draw (0,0) node[above right]{O}; 
        \draw[very thick,domain=0:4.5] plot(\x,0.5*\x*\x-0.5*\x*3-0.875) node[left]{Left-hand side of \Cref{eq:QED}};
        \draw[very thick, dashed, domain=-1.5:0] plot(\x,0.5*\x*\x-0.5*\x*3-0.875);
        \draw[very thick,domain=0:4.5] plot(\x,0.88) node[above left, shift={(-1,0)}]{Right-hand side of \Cref{eq:QED}};
        \draw[very thick, domain=-1.5:0] plot(\x,0.88) ;
        \draw[thick] (1.5,-0.1)--(1.5,0.1) node[above]{$\frac{-\tau_1 \mathbb{E}\lrs{u}}{2\tau_0}$};
        \draw[thick] (3.9,0.1)--(3.9,-0.1) node[below]{$n_\mathrm{cr}$};
        \draw[thick, dashed] (3.9,0.88)--(3.9,0); 
        \draw[thick] (0.1,-0.875)--(-0.1,-0.875) node[left]{$\tau_2 \mathbb{E}\lrs{u^2}<0$};
        \draw[thick, dashed] (1.5,-2)--(1.5,0.); 
    \end{tikzpicture}
    \caption{Schematic view of the unique real solution of \Cref{eq:QED}}
    \label{fig:schematicview2}
\end{figure}

This section demonstrates the uniqueness of the surf-riding threshold in the framework of Melnikov's method.
First, the authors start with the quadratic $K_\mathrm{T}\lrp{J}$, which is adopted in the level-2 vulnerability criterion for broaching in the SGISC.
Substituting \Cref{eq:Te} into \Cref{eq:melbalance2} yields
\begin{equation} \label{eq:MelQ}
    \frac{1}{2 \pi} \int_{-\pi}^\pi \lrp{\tau_0 n_\mathrm{cr}{}^2 + \tau_1 u n_\mathrm{cr} +\tau_2 u^2}\dif{y} = \frac{1}{2 \pi} \int_{-\pi}^\pi R\lrp{u}\dif{y}
\end{equation}
The left-hand side is the quadratic function with $n_\mathrm{cr}$.
Since the following definite integrals with respect to $y$ are calculated as
\begin{align}
    \frac{1}{2 \pi} \int_{-\pi}^\pi \dif{y} &= 1 \\
    \frac{1}{2 \pi} \int_{-\pi}^\pi u \dif{y} &= 
    c_\mathrm{w}-\frac{4}{\pi} \sqrt{\frac{f}{k\lrp{m+m_\mathrm{x}}}} \triangleq \mathbb{E}\lrs{u} \\
    \frac{1}{2 \pi} \int_{-\pi}^\pi u^2 \dif{y} &= c_\mathrm{w}{}^2 -\frac{8}{\pi} c_\mathrm{w} \sqrt{\frac{f}{k\lrp{m+m_\mathrm{x}}}} + 2 \frac{f}{k\lrp{m+m_\mathrm{x}}} \notag \\
    & \triangleq \mathbb{E}\lrs{u^2} > 0
\end{align}
\Cref{eq:MelQ} is rewritten as
\begin{align} 
    & \tau_0 n_\mathrm{cr}{}^2 + \tau_1 \mathbb{E}\lrs{u} n_\mathrm{cr}
    + \tau_2 \mathbb{E}\lrs{u^2} =\frac{1}{2 \pi} \int_{-\pi}^\pi R\lrp{u}\dif{y} \notag \\
    \Leftrightarrow ~ & \tau_0 
            \lrp{ n_\mathrm{cr} - \frac{- \tau_1 \mathbb{E}\lrs{u}}{2\tau_0}}^2 - \tau_0 \lrc{ \frac{\lrp{\tau_1 \mathbb{E}\lrs{u}}^2}{4\tau_0{}^2} + \frac{-\tau_2 \mathbb{E}\lrs{u^2}}{\tau_0} 
        }\notag\\
    & \phantom{\tau_0 n_\mathrm{cr}{}^2 + \tau_1 \mathbb{E}\lrs{u} n_\mathrm{cr}
    + \tau_2 \mathbb{E}\lrs{u^2}}
    =\frac{1}{2 \pi} \int_{-\pi}^\pi R\lrp{u} \dif{y} \label{eq:QED}
\end{align}
Since $\tau_0$ is positive as shown in \Cref{eq:tau0}, and since $\tau_2$ is negative as shown in \Cref{eq:tau2}, the left-hand side is a convex downward quadratic function that is negative when $n_\mathrm{cr}$ is zero.
The schematic view of the solution of \Cref{eq:QED} is shown in \Cref{fig:schematicview2} as the crossing points of the quadratic function (the left-hand side of \Cref{eq:QED}) and the straight line (the right-hand side of \Cref{eq:QED}). 
Obviously, \Cref{eq:QED} has two distinct real roots, one positive and the other negative if 
\begin{equation} \label{eq:rtcond}
    \frac{1}{2 \pi} \int_{-\pi}^\pi R\lrp{u}\dif{y} > \tau_2 \mathbb{E}\lrs{u^2}
\end{equation}
However, when $n_\mathrm{cr}$ is negative and $u$ is positive, that is, adverse rotation of the propeller with the ship moving forward, the effective thrust of the propellers becomes negative and the $K_\mathrm{T}$ for $n \in \mathbb{R}_{\geq 0}$ is ineffective.
Therefore, \Cref{eq:QED} always has a unique real $n_\mathrm{cr}$ regardless of the value of $\tau_1$.

In the above discussion, it was demonstrated that the surf-riding threshold can be uniquely determined if the $K_\mathrm{T}$ is approximated under the condition of $n > 0$ and \Cref{eq:rtcond}. 
Further discussion of the equivalent wave condition under which the unique surf-riding threshold is given is provided in the appendix, and it is confirmed that \Cref{eq:rtcond} is satisfied in most waves.

\section{Physical Interpretation of the Uniquness of the surf-riding threshold}\label{sec:threshold}
\begin{figure}[h]
    \centering
    \begin{tikzpicture}
        \draw[->,>=stealth,very thick] (0,0)--(6,0) node[right]{$n$}; 
        \draw[->,>=stealth,very thick] (0,-2.5)--(0,3); 
        \draw (0,0) node[above right]{O}; 
        \draw[very thick,domain=0:4.9] plot(\x,0.2*\x*\x+0.04*\x-2) node[below left, shift={(-0.15,0)}]{$\frac{1}{2\pi} \int_{-\pi}^\pi T_\mathrm{e}\lrp{u;n} \dif{y}$};
        \draw[very thick,domain=0:6] plot(\x,1.5) node[above left, shift={(-3,0)}]{$\frac{1}{2\pi} \int_{-\pi}^\pi R\lrp{u} \dif{y}$};
        \draw[thick, dashed] (4.0845,1.5)--(4.0845,0); 
        \draw[thick] (4.0845,0.1)--(4.0845,-0.1) node[below]{$n_\mathrm{cr}$};
        \draw[thick] (0,-2) node[below right]{$\frac{1}{2\pi} \int_{-\pi}^\pi \left. T_\mathrm{e}\lrp{u;n} \right|_{n=0} \dif{y}<0$};
    \end{tikzpicture}
    \caption{Schematic view of the unique real solution of \Cref{eq:melbalance2}}
    \label{fig:schematicview3}
\end{figure}
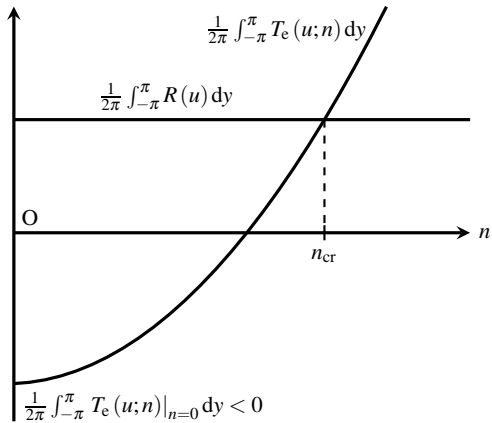

The relationship between $T_\mathrm{e}$ and $R\lrp{u}$, as stated in \Cref{eq:melbalance2}, is used to discuss the physical interpretation of the unique surf-riding threshold. 
In the context of surf-riding, it is reasonable to assume $u,n \in \mathbb{R}_{\ge 0}$. 
Under these conditions, $T_\mathrm{e}$ can be assumed to have the following properties.
\begin{enumerate}
    \item $T_\mathrm{e}$ increases monotonically with $n$. \label{enu:Te}
    \item When $n=0$, $T_\mathrm{e}$ has a negative finite value. \label{enu:Te0}
\end{enumerate}
As in \Cref{sec:unique}, the definite integral for the resistance on the right-hand side is constant with $n$. 
From the first point, the definite integral on the left-hand side of \Cref{eq:melbalance2} increases monotonically with respect to $n$.
The second point corresponds to the fact that the propeller generates a negative finite thrust when $u>0$ and $n = 0$. 
Therefore, \Cref{eq:melbalance2} always has a unique positive solution. 
The schematic view is shown in \Cref{fig:schematicview3}.
Although $\frac{1}{2\pi} \int_{-\pi}^\pi T_\mathrm{e}\lrp{u;n} \dif{y}$ does not always increase monotonically with small $n$ in \Cref{fig:schematicview2}, this is simply due to the quadratic fitting of $K_\mathrm{T}$ with respect to $J$.

\section{Conclusions}\label{sec:conclusion}
    The authors discussed the uniqueness of the approximate analytical solution of the surf-riding threshold in the SGISC.
    The authors interpreted Melnikov's method in physical terms, viewing it as a balance between the work done by effective thrust and resistance, and found that 
    \begin{enumerate}
        \item the surf-riding threshold $n_\mathrm{cr}$ is obtained as the solution of a quadratic equation when $K_\mathrm{T}$ is expressed as a quadratic function of $J$ as adopted in SGISC, which has two distinct real roots, one positive and the other negative;
        \item the positive root is the unique surf-riding threshold if the right-hand side of \Cref{eq:QED} is positive; 
        \item in most of waves, the right-hand side of \Cref{eq:QED} is positive; and
        \item the above does not depend on the approximation of the calm-water resistance $R$.
    \end{enumerate}
    
\begin{acknowledgements}
    The authors express great appreciation to Prof.~Atsuo Maki of Osaka Univerity for constructive suggestions and discussions.
    This study was supported by a Grant-in-Aid for Scientific Research from the Japan Society for the Promotion of Science (JSPS KAKENHI Grant JP21K14362). 
\end{acknowledgements}

%
\section*{Conflict of interest}
    The authors declare that they have no conflict of interest.

\bibliographystyle{spbasic}      
\bibliography{main}   

\section*{Appendix} \label{sec:app}
The unique surf-riding threshold is given when \Cref{eq:rtcond} is satisfied
under \Cref{eq:uint} in the framework of Melnikov's method.
Since the right-hand side of \Cref{eq:rtcond} is negative, 
\begin{equation} \label{eq:cond}
     \int_{-\pi}^\pi R\lrp{u}\dif{y} > 0
\end{equation}
can be considered as the necesarry condition of \Cref{eq:rtcond}.
It is obvious that \Cref{eq:cond} is satisfied when $u>0$ along the heteroclinic orbit. Therefore, from \Cref{eq:uint}, 
\begin{equation} \label{eq:uconst}
	f < \frac{1}{4} \lrp{m + m_\mathrm{x}} g
\end{equation}
Assuming that the calm-water resistance has an equal magnitude for $u<0$, \Cref{eq:cond} is almost equivalent to
\begin{equation} \label{eq:uconst2}
	\mathbb{E}\lrs{u}>0 \quad \Leftrightarrow \quad f < \frac{\pi^2}{16} \lrp{m + m_\mathrm{x}} g
\end{equation}
In reality, the calm-water resistance would be larger for $u<0$, so the condition for $f_\mathrm{R0}$ where 
\begin{equation} \label{eq:r0}
	\int_{-\pi}^\pi R\lrp{u}\dif{y} = 0
\end{equation}
is satisfied ranges
\begin{equation}
	\frac{1}{4} \lrp{m + m_\mathrm{x}} g < f_\mathrm{R0} < \frac{\pi^2}{16} \lrp{m + m_\mathrm{x}} g
\end{equation}

The wave-induced surge force is calculated as the Froude-Krylov force in the SGISC. Therefore, the amplitude is calculated as
\begin{align}
    f &= \mu \cdot \rho g k \frac{H}{2} \sqrt{F_\mathrm{c}{}^2+F_\mathrm{s}{}^2} \notag \\
      &= \mu \cdot \pi \rho g \frac{H}{\lambda} \sqrt{F_\mathrm{c}{}^2+F_\mathrm{s}{}^2}
\end{align}
where
\begin{align}
    F_\mathrm{c} &= \int S\lrp{x} \sin \lrp{kx} e^{-k\frac{d\lrp{x}}{2}} \dif{x} \\
    F_\mathrm{s} &= \int S\lrp{x} \cos \lrp{kx} e^{-k\frac{d\lrp{x}}{2}} \dif{x} \\
    m            &= \rho \int S\lrp{x} \dif{x}
\end{align}
Here, $\mu$ is the correction factor due to the diffraction effect, $H$ is the wave height, $x$ is the body-fixed coordinate system along the transverse direction ($x=0$ indicates the midship), $S$ is the area of the submerged portion at $x$ and $d$ is the draught at $x$.
In the vulnerability criterion of the SGISC, the diffraction effect on the wave-induced surge force is neglected and is simply considered as $\mu = 1$. 
For operational guidance, the explanatory note to the SGISC suggests numerical estimation of $\mu$ or the use of experimental correction factor \cite{Ito2014-hb,Imo2022-uf} as
\begin{equation}
    \mu = \left\{ 
    \begin{array}{l@{\quad}l}
        1.46 C_\mathrm{b} - 0.05 & \lrp{C_\mathrm{m}<0.86} \\
        \lrp{5.76-5.00 C_\mathrm{m}} C_\mathrm{b} - 0.05 & \lrp{0.86 \geq C_\mathrm{m} \geq 0.94} \\
        1.06 C_\mathrm{b} - 0.05 & \lrp{C_\mathrm{m} > 0.94}  
    \end{array} 
    \right. 
\end{equation}
In simple terms, the value is almost $\mu = 0.7$ for the operational guidance.

An upper bound on $f$ is considered with maximised wave slope acting on each submerged section as
\begin{align} \label{eq:fandS}
    f &< \mu \cdot \pi \rho g \frac{H}{\lambda} \int S\lrp{x} \cdot  1 \cdot 1 \dif{x} \notag \\
      &= \mu \cdot \pi \frac{H}{\lambda} m g
\end{align}
Since $m_\mathrm{x} \simeq 0.1 m$, \Cref{eq:fandS} can be rewritten as
\begin{equation}
    f < \mu \cdot \frac{4\pi}{1.1} \frac{H}{\lambda} \cdot \frac{1}{4} \lrp{m + m_\mathrm{x}} g
\end{equation}
\begin{equation}
    f < \mu \cdot \frac{16}{1.1 \pi} \frac{H}{\lambda} \cdot \frac{\pi^2}{16} \lrp{m + m_\mathrm{x}} g
\end{equation}
As a result, \Cref{eq:uconst,eq:uconst2} is almost equivalent to the following wave condition
\begin{align}
    \mu \cdot & \frac{4\pi}{1.1} \frac{H}{\lambda} < 1 \notag \\
              &  \Rightarrow  \left\{ 
              \begin{array}{l@{\quad}l}
                  \frac{H}{\lambda} < \frac{1}{11.4} & \lrp{\mu = 1 \text{, vulnerability criterion}} \\
                  \frac{H}{\lambda} < \frac{1}{8.00} & \lrp{\mu = 0.7 \text{, operational guidance}} 
              \end{array} 
              \right. \\
    \mu \cdot & \frac{16}{1.1\pi} \frac{H}{\lambda} < 1 \notag \\
              &  \Rightarrow  \left\{ 
              \begin{array}{l@{\quad}l}
                  \frac{H}{\lambda} < \frac{1}{4.63} & \lrp{\mu = 1 \text{, vulnerability criterion}} \\
                  \frac{H}{\lambda} < \frac{1}{3.24} & \lrp{\mu = 0.7 \text{, operational guidance}} 
              \end{array} 
              \right.
\end{align}
Therefore, the wave condition for $\lrp{\frac{H}{\lambda}}_\mathrm{R0}$ where \Cref{eq:r0} is satisfied is
\begin{equation} 
	\left\{ 
          \begin{array}{l@{\quad}l}
              \frac{1}{11.4} < \lrp{\frac{H}{\lambda}}_\mathrm{R0} < \frac{1}{4.63} & \lrp{\mu = 1 \text{, vulnerability criterion}} \\
              \frac{1}{8.00} < \lrp{\frac{H}{\lambda}}_\mathrm{R0} < \frac{1}{3.24} & \lrp{\mu = 0.7 \text{, operational guidance}} 
          \end{array} 
          \right. 
\end{equation}
In the vulnerability criterion, the maximum wave steepness is $0.15 \simeq 1/6.67$. 
Therefore, there is the possibility that $u < 0$. 
However, \Cref{eq:rtcond} is usually satisfied.

\end{document}